# Functional Data Analysis in Electronic Commerce Research

**Wolfgang Jank and Galit Shmueli**

*Abstract.* This paper describes opportunities and challenges of using functional data analysis (FDA) for the exploration and analysis of data originating from electronic commerce (eCommerce). We discuss the special data structures that arise in the online environment and why FDA is a natural approach for representing and analyzing such data. The paper reviews several FDA methods and motivates their usefulness in eCommerce research by providing a glimpse into new domain insights that they allow. We argue that the wedding of eCommerce with FDA leads to innovations both in statistical methodology, due to the challenges and complications that arise in eCommerce data, and in online research, by being able to ask (and subsequently answer) new research questions that classical statistical methods are not able to address, and also by expanding on research questions beyond the ones traditionally asked in the offline environment. We describe several applications originating from online transactions which are new to the statistics literature, and point out statistical challenges accompanied by some solutions. We also discuss some promising future directions for joint research efforts between researchers in eCommerce and statistics.

*Wolfgang Jank is Assistant Professor of Management Science and Statistics, Department of Decision and Information Technologies, Robert H. Smith School of Business, University of Maryland, College Park, Maryland 20742, USA e-mail: wjank@rhsmith.umd.edu. Galit Shmueli is Assistant Professor of Management Science and Statistics, Department of Decision and Information Technologies, Robert H. Smith School of Business, University of Maryland, College Park, Maryland 20742, USA e-mail: gshmueli@rhsmith.umd.edu.*

## 1. INTRODUCTION

Functional data analysis (FDA) has been gaining momentum in many fields. While much of the methodological advances have been made within the statistics literature, FDA has found many useful applications in the agricultural sciences (Ogden et al., 2002), the behavioral sciences (Rossi, Wang and Ramsay, 2002), in medical research (Pfeiffer et al., 2002) and many more. One reason for this momentum is the technological advancement in computer storage and computing power. Today's researchers gather more and more data, often automatically, and store them in large databases. However, these new capabilities for data generation and data storage have also led to new data structures which do not necessarily fit into the classical statistical concept. Researchers measure characteristics of customers over time, store digitalized two- or three-dimensional images of the brain, and record three- or even four-dimensional movements of objects through space and time. Many of these new data structures call for new statistical methods in order to unveil the information that they carry. Data can contain trends that vary in longitudinal or spatial aspects, that vary across different groups of customers or objects, or that show different magnitudes of dynamics.

FDA is a tool-set that, although based on the ideas of classical statistics, differs from it (and, in a sense, generalizes it), especially with respect to the type of data structures that it encompasses. While the underlying ideas for FDA have been around for a longer time, the surge in associated research can be attributed to the monograph of Ramsay and Silverman ([1997](#)). In FDA, the object of interest is a set of curves, shapes, objects, or, more generally, a set of *functional observations*. This is in contrast to classical statistics where the interest centers around a set of data vectors. In recent years, a range of classical statistical methods have been generalized to the functional framework; James, Hastie and Sugar ([2000](#)) developed a principal components approach for a set of sparsely sampled curves. Other exploratory tools include curve clustering (see Abraham et al., [2003](#); James and Sugar, [2003](#); Tarpey and Kinateder, [2003](#)) and curve classification (see Hall, Poskitt and Presnell, [2001](#); James and Hastie, [2001](#)). Classical linear models have also been generalized to functional ANOVA (Fan and Lin, [1998](#); Guo, [2002](#)), functional regression (Faraway, [1997](#); Cuevas, Febrero and Fraiman, [2002](#); Ratcliffe, Leader and Heller, [2002](#)) and the functional generalized linear model (Ratcliffe, Heller and Leader, [2002](#); James, [2002](#)). Moreover, Ramsay ([2000b](#)) and Ramsay and Ramsay ([2002](#)) suggest differential equations for data of functional form. While this list is far from complete, it shows some of the current methodological efforts in this emerging field.

Electronic commerce (eCommerce) is a growing field of scholarly research especially in information systems, economics and marketing, but it has received little to no attention in statistics. This is surprising because it arrives with an enormous amount of data and data-related questions and problems. Like other web-based data, eCommerce data tend to be very rich, clean and structurally different from offline data. eCommerce research arrives with many new data- and model-related challenges that promise new ideas and motivation for further methodological advancements of FDA. One of the main characteristics of eCommerce data is the combination of longitudinal information (time-series data) with cross-sectional information (attribute data). A sample of $n$ records typically comprises $n$ time series, each linked with a set of $n$ attributes. Take *eBay*'s online auctions as an example. There, each auction is characterized by a time series of bids placed over time. This information is coupled with additional auction attributes such as a seller's rating, the auction duration and the currency used. Another example is online product ratings on Amazon.com or movie ratings on *Yahoo! Movies* (see Dellarocas and Narayan, [2006](#)). Yahoo! Movies allows users to rate any movie according to different measures. This results in a time series that describes the average daily rating or the number of daily postings (or both) from the date of the movie release until the time of data collection. This information is coupled with attribute data about the movie such as the movie genre and critics' rating. A third example, described in Stewart, Darcy and Daniel ([2006](#)), is the evolution of open-source software projects that is monitored by websites such as [SourceForge.net](#). Here an observation is a certain project, and it is characterized by a time series that describes project complexity from its first release until the time of data collection. Each project also has associated attributes such as the number of developers, the operating system used and the programming language.

The combination of longitudinal and cross-sectional information is only one typical aspect of eCommerce data. Another aspect is the uneven spacing between events. In many cases, the observed time series is composed of events influenced by multiple users or agents who access the web at different points in time (and from different geographical locations). Consequently, the resulting times when new events arrive are extremely unevenly spaced. This is in contrast to traditional time series, which are typically recorded at predefined and equidistant time-points, such as daily, monthly or quarterly scales. Furthermore, because of psychological, economic or other reasons, eCommerce time series tend to feature very sparse areas at some times, followed by extremely dense areas at other times. For instance, bidding in eBay auctions tends to be concentrated at the end, resulting in very sparse bid-arrivals during most of the auction except for its final moments, where the bidding volume can be extremely high.

eCommerce not only creates new data challenges, it also motivates the need for innovative models.

While the field of economics has created many theories for understanding economic behavior at the individual and market level, many of these theories were developed before the emergence of the World Wide Web. The existence of the web now allows researchers, for the first time, to observe and record data about economic behavior on a large-scale basis. As it turns out, however, observed data often do not support classical economic theories. As a result, empirical research is thriving. In fact, the empirical literature has continuously shown that online behavior deviates in many ways from offline behavior and from what is expected by economic theory. This calls for new economic models that can be validated empirically. In addition, the availability of eCommerce data allows researchers to ask new types of questions. One major enhancement is the ability to study not only the evolution of a process, but also its dynamics: how fast it moves and how suddenly it changes, its rate of change and how this rate differs at different time-points. Studying dynamics of processes can be very relevant in the online world, because it allows new approaches for characterizing eCommerce processes (and thus distinguishing between diverse processes), and even forecasting them (Wang, Jank and Shmueli, 2006). Changing dynamics are inherent in a fast-moving environment like the online world. Fast movements and change imply nonstationarity which poses challenges to traditional time series modeling. And finally, it is important to point out that for any one process that we observe in the online world, there typically exist many, many replicates of the same (or at least very similar) process. On eBay, for instance, if we think of the formation of price between the start and the end of an auction as a process of interest, then there exist several million similar processes of that form, taking place at any given day on eBay. The replication of processes, or time series, fits naturally within the FDA framework and makes this an ideal ground for the advancement of new functional methodology.

Finally, eCommerce typically arrives with huge databases which can put a computational burden on users' storage and processing facilities. This burden is often increased by the complicated structure of eCommerce data. Taking a functional data approach, one can relieve some of that burden. FDA operates on functional objects which can be more compactly represented than the original data. Taking a functional approach may therefore be advantageous also from a resource point of view.

The process of studying a set of data via functional methods consists of two principal steps: First, the functional object is "recovered," typically by means of smoothing. There are multiple different ways in which this smoothing step can be executed, and there are many challenges during that step. Second, the resulting functional object is used for data exploration and analysis. Exploratory data analysis (including data visualization and summary) is performed in order to learn about general characteristics as well as unusual features and anomalies in the data. Analysis includes explanatory and predictive modeling and inference, just as in classical statistics. In the next sections we focus on the challenges and problems that arise during these steps within the eCommerce context. We would like to note that our point of view of the functional approach and its application to eCommerce has been forged during the teaching of so-called *Research Interaction Teams* (www.amsc.umd.edu/Courses/RITDescrips/HowAndWhy.html) which are research classes that involve graduate students from the Statistics and the Applied Mathematics and Scientific Computation programs at the University of Maryland. Several of our studies performed during these classes have led to new methodological and practical insights.

## 2. RECOVERING FUNCTIONAL OBJECTS

The first step in any functional data analysis consists of recovering, from the observed data, the underlying functional object. There exist a variety of methods for recovering functional objects from a set of data, all of which are typically based on some kind of smoothing. As a result of the smoothing, and of characterizing the smooth object by its smoothing parameters only, we obtain a low-dimensional functional object. We focus here on objects, and in particular curves, that are based on unevenly spaced time series and of which we have multiple replications. An example is a set of bid histories from eBay auctions, as shown in Figure 1. The four panels correspond to four separate seven-day auctions for a new Palm PDA. Each consists of the bids (in $) placed at different times during the auction.

### 2.1 Challenges in Choosing the Right Smoother

The first step in recovering the functional object is to choose a family of basis functions. The choice of the basis function depends on the nature of the data, on the level of smoothness that the application warrants, on what aspects of the data we want to study, on the size of the data and on the types of analyses that we plan to perform. For example, to represent the price path of an online auction for the purpose of, say, studying price dynamics, one could use monotone smoothing splines (Ramsay, 1998) since prices in auctions increase monotonically. Besides maintaining the price monotonicity, this approach also permits the computation of derivatives which lend themselves to price dynamics. However, fitting monotone splines is computationally more intensive than fitting ordinary polynomial smoothing splines. In that sense, it may prove impractical to compute monotone splines for very large databases if time and memory restrictions exist. In addition, polynomial smoothing splines can be represented as a linear combination of basis functions. The practical meaning of this is that if we use polynomial smoothing splines and the intended analysis is based on a linear operation (such as computing average curves, fitting functional linear regression models, or performing functional principal components analysis), we can operate directly on the basis function coefficients without any loss of information. The same operation using monotone splines would require an approximation step due to the need to first represent the continuous curve in a finite-dimensional manner by evaluating it on a grid. Conversely, if the type of operation is nonlinear, then one would have to perform a grid-based computation for either type of spline and the choice would therefore not matter from this point of view. Thus, the way we recover the functional object is strongly influenced by a variety of different objectives all of which might compete with one another.

Recovering functional objects often involves more than deciding on the appropriate type of smoother. This can include a preprocessing step via interpolation, thereby creating a *raw functional* (e.g., Ramsay and Silverman, 2002, page 21). This alleviates the problem of the unevenly spaced series that are common in eCommerce. An important aspect in any functional data analysis is the robustness of *analysis results* with respect to the choice and level of smoothing. A general study of this sort was carried out by our research interaction team, comparing the effects of smoothing splines versus monotone splines on the conclusions derived from a functional regression on the price path in online auctions (the functional object) as a function of explanatory variables such as the seller rating and the opening bid (both scalar) and current number of bids (a functional explanatory variable). The study indicates that both smoothers lead to similar conclusions (Alford and Urimi, 2004).

Another example of the challenges in choosing an appropriate smoothing method is the functional representation of online movie ratings. By that we mean the series of user movie ratings on online services such as Yahoo.com. The volume of user postings is highly periodic, with heavier activity on weekends (when people tend to watch movies in the theaters). Fourier basis functions were found to be a better choice among different alternatives for capturing the cyclical posting patterns (Wu, 2005).

## 2.2 Additional Data Challenges

Our different studies using eCommerce data have raised further challenges in the functional object recovery stage that have previously not been addressed in the literature. The first such challenge is handling the extremely unevenly distributed measurements in eCommerce data. That is, the number and location of events vary drastically from one functional object to another. One typical example is the bid arrival in eBay auctions. Returning to Figure 1, it can be seen that some bid histories are very dense at the auction end, while others are much sparser, and in addition the overall number of bids per auction can vary widely between none in some auctions, and more than 100 in others. And yet, while the varying number of bids per auction may suggest the use of a varying set of smoothing methods, we prefer the use of a *single* family of smoothers. The reason for this is that, in the end, the choice of the smoother is merely a means to the end of arriving at a unifying functional object and it is not the direct object of our interest. Coming back to the example of online auctions with sometimes few and sometimes many bids, this motivates the need for new methodological advances in creating functional objects that naturally incorporate all extremes under one hat. Some promising approaches in that direction can be found in James and Sugar (2003) and James, Hastie and Sugar (2000).

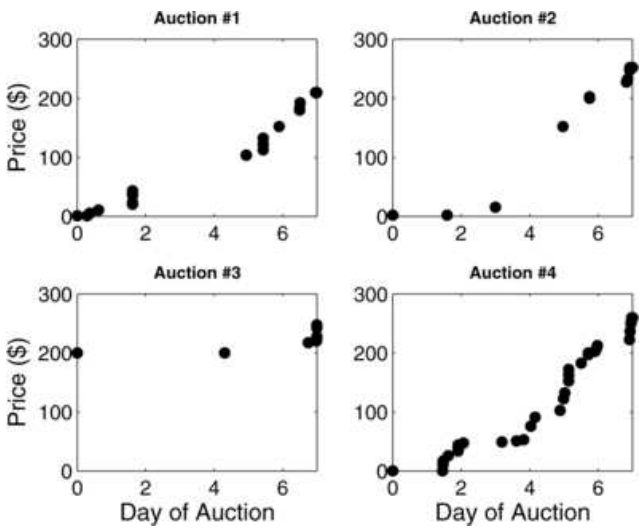

Fig. 1. *Scatterplots describing the bid history in each of four eBay auctions, each lasting seven days.*

The extreme structure of eCommerce data is challenging even for very basic visualization tasks: standard time-series visualization tools typically require evenly spaced events! Because of this restrictive requirement, we collaborated with colleagues at the *Human–Computer Interaction Lab* at the University of Maryland to develop new interactive visualization tools that can accommodate this special data structure. We evaluated different approaches for representing eBay bid histories by an evenly spaced equivalent without losing important information with respect to order, magnitude and distance between the bids (see Aris et al., 2005). Interestingly, the final choice was to use a functional approach by first smoothing the bid history, and then feeding an evenly spaced grid of the smooth curves (and their derivatives) into a standard visualization tool.

Another challenge typical to eCommerce data is defining meaningful start and end points of the functional objects in order to align the curves. The problem of aligning functional objects is related to the problem of registration (Ramsay and Silverman, 2005, Chapter 7), but there are several additional complications here: Many web-based events do not start and end at the same time. For instance, online product ratings over time have different starting points, depending on when the product was first released to the market, when the first rating was placed, etc. They also often have different ending points, for instance, if one product is prematurely taken off the market, if it is replaced by a product-upgrade, and so on. Another issue that complicates object alignment is that the data collection process itself may act as a censoring mechanism. Therefore, it is not obvious how methods such as *landmark registration*, where curves are aligned according to one particular feature of the curve such as its peak, can be adapted to handle this situation. Finally, selecting the units for the time axis can be challenging. In some applications calendar time (e.g., the date and time a transaction took place) is reasonable, whereas in other applications the event index (i.e., the order of the event arrival) might make more sense. And yet in other cases an entirely different "clock" would be even more suitable. For instance, we pointed out that eBay auctions typically exhibit very low bidding activity during most of the auction and then extremely high activity near the end. For this reason an auction might be better represented by a clock that "shrinks" the low-activity period and "stretches" the high-activity period, thereby putting more emphasis on the part that matters more. All these issues are illustrated and discussed further in the paper by Stewart, Darcy and Daniel (2006).

## 3. FUNCTIONAL EDA

After the data are represented by functional objects, the analysis steps follow the same process as in classical statistics, with the first step being exploratory data analysis (EDA). EDA includes data summaries, visualization, dimension reduction, outlier detection, and more. The main difference between FDA and classical statistics is the way in which the methods are applied and especially how they are interpreted.

### 3.1 Static vs. Interactive Visualization

Starting with visualization, our goal is to:

1. Visualize a sample of curves.
2. Inspect summaries of these curves.
3. Explore conditional curves, using various relevant predictor variables.

To that end, one solution is to create *static* graphs. For instance, Figure 2 shows the price evolution in 34 eBay auctions for various magazines. We can see large heterogeneity across the price formation process at different times of the auction. We refer to this approach as static, since once the graph is generated it can no longer be modified by the user without running the software code again. This static approach is useful for differentiating subsets of curves by attributes (e.g., by using color), or for spotting outliers. However, a static approach does not allow for an interactive exploration of the data. By interactive we mean that the user can perform operations such as zooming in and out, filtering the data and obtaining details for the filtered data, and do all of this from within the graphical interface. Interactive visualizations for the special structure of eCommerce data are not straightforward, and solutions have only been proposed recently (Aris et al., 2005; Shmueli et al., 2006). One such solution is *Auction-Explorer* (www.cs.umd.edu/hcil/timesearcher), which is tailored to handle the special structure of online auction data. A snapshot of its user interface is shown in Figure 3. The interface includes several panels, which correspond to the price curves (top left), their dynamics (not shown in this view) and the corresponding attribute data (top right). The curves can be filtered to display subsets according to a selection of attribute values, according to a selection of curves, and one can do pattern search. Summarization is achieved through on-the-fly summary statistics for attributes, and a "streaming boxplot" called a "river plot" of the curves (bottom panel of Figure 3). This is yet another attempt to generalize classical visualization methods to the functional environment.

### 3.2 Data Reduction

Another goal of EDA is data reduction. Two of the methods that are useful in this context are curve clustering and functional principal components analysis. Curve clustering partitions the set of curves into a few clusters, thereby reducing the space of observations, and attempts to derive insight from the resulting clusters. The clustering can be applied to the curves themselves or to their derivatives. Jank and Shmueli (2005) apply curve clustering to bid histories of eBay auctions and find three main clusters. Linking the curve information with attribute information, they find that the different clusters correspond to three types of auctions: "*greedy sellers*,"

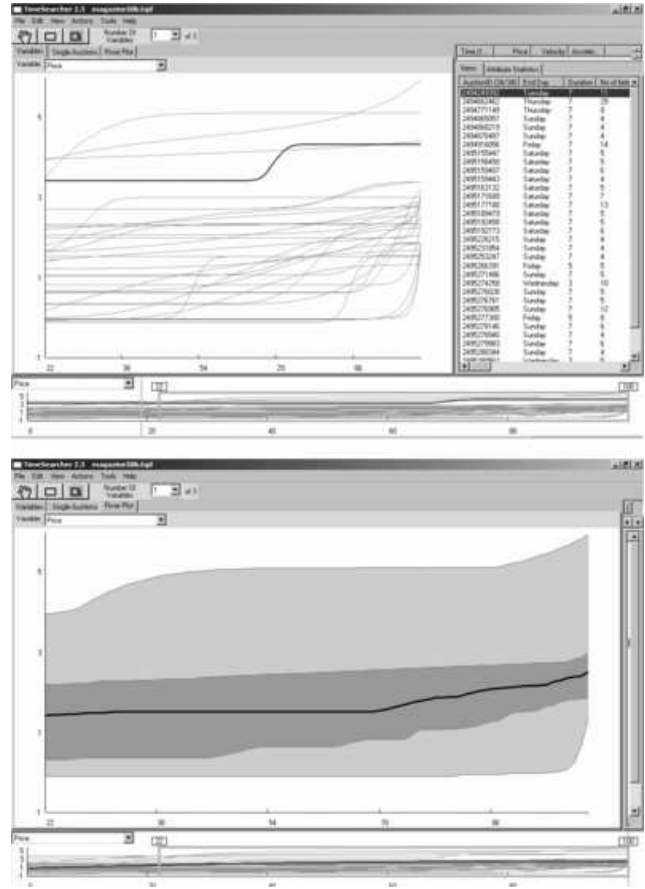

"*bazaar auctions*" and "*experienced seller/buyer*" auctions. Each of these types characterizes a different auction profile, combining static and dynamic information. For instance, "*greedy seller*" auctions have the highest average opening price and the lowest closing price. Unsurprisingly, they do not attract much competition, since unjustified high opening prices tend to deter users from bidding. Low competition is also known to lead to lower prices. Sellers in these auctions are, on average, less experienced than those in other clusters (as can be measured by their eBay rating), scheduling most auctions to end on a weekday. These auctions also attract experienced winners who take advantage of the poorer auction design and resulting lower prices. The price dynamics of this cluster reflect this setting: price starts accelerating late in the auction, not allowing it to achieve its full impact by the time the auction closes. This mix of insight into static seller and bidder characteristics coupled with the price dynamics is only available with a functional approach.

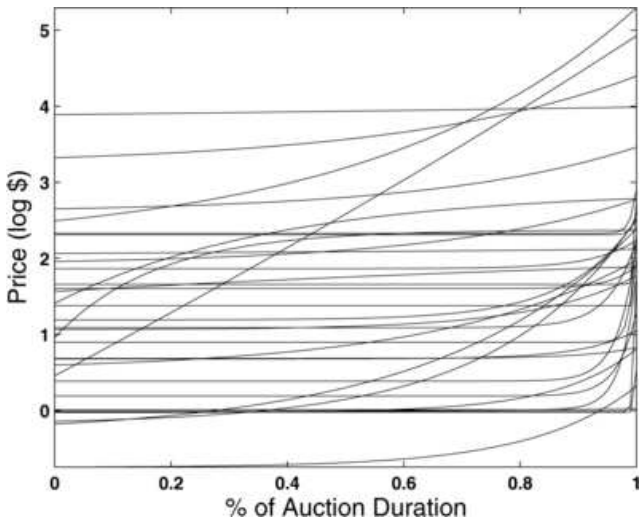

Fig. 2. *Static plot of the price progression in 34 eBay auctions for magazines.*

Another popular method is functional principal components analysis (f-PCA). The method uses standard PCA to find principal sources of variability in curves (or other functional objects). If we consider curves that represent a process over time, then f-PCA can help us find "within-curve" (or more generally, "within-process") variation, thereby condensing the time axis. This is done by selecting a discrete grid of time-points and treating the points as the variables in ordinary PCA. A preliminary study by Hyde, Moore and Hodge (2004) applied f-PCA to price curves and derivative curves of a sample of eBay auctions for premium wristwatches. They found that two or three principal components captured most of the within-curve variation: One source is price variation during the middle of the auction and the other distinguishes the price variability between the beginning and end of the auction. Similar results were obtained when using the price-dynamics curves. Hyde, Moore and Hodge (2004) also used f-PCA to compare sources of "within-process" variation across different brands for the same product category as well as for different product categories. It was found that price is most uncertain during mid-auction. As the auction approaches its end, though, the price becomes more predictable, especially in common-value auctions (see also Wang, Jank and Shmueli, 2006).

An alternative approach of principal components analysis to functional data that, to the best of our knowledge, has not been explored would be to treat the observations as the dimension to be transformed. The idea is to find main sources of variation across curves (instead of within curves), achieving a goal similar to curve clustering, where main features of the curves are highlighted. The exact meaning and interpretation of this variation deserve further attention.

## 4. FUNCTIONAL MODELING, INFERENCE AND PREDICTION

There is quite a lot of ongoing research on generalizing classical regression models to the functional setting. Examples include linear regression with functional predictors (Ratcliffe, Leader and Heller, 2002) or a functional response (Faraway, 1997), logistic regression (Ratcliffe, Heller and Leader, 2002), functional linear discriminant analysis (James and Hastie, 2001) and general linear models with functional predictors (James, 2002).

### 4.1 Information in eCommerce Processes

Current empirical research in eCommerce relies on the use of very standard statistical tools such as least-squares regression. These tools are used to investigate how, say, the closing price in an online auction relates to other auction-specific information. To that end, one sets the closing price as the response variable, and regresses it on potential explanatory variables such as the opening bid, the auction duration, the seller rating, etc. (see, e.g., Bajari and Hortaçsu, 2003, 2004; Lucking-Reiley et al., 2000). While this approach is certainly useful for understanding *some* of the variation in closing prices, it also leads to loss of a large amount of potentially useful information about everything that happened between the start and end of the auction. More generally, current research uses a response variable that is an aggregation of the process of interest: the maximum bid in online auctions, the average product rating, etc. This (direct or indirect) choice is guided by economic importance but is also likely done so that standard models can be applied. Furthermore, the choice of independent variables is limited to static "snapshot" information. The existence of more detailed data, however, can potentially shed more light on the entire process rather than only its aggregated form.

In the online auction example, variables like the opening bid, the auction duration and the seller's rating are determined *before* the auction start and thus do not capture any of the information that arrives after. However, it is well-known that events that occur *during* the auction can also affect the final price. For instance, the number of competing bidders, the bidders' experience and the bid timings can influence the final price. These three variables are available only *after* the auction starts and in fact the information they carry changes as the auction progresses.

While it is possible to include time-varying explanatory variables like the number of bidders into a regression model, such a model would no longer be considered "standard" in the classical least-squares sense since it would have to account for time-dependence between the explanatory variable and the response, and also within the explanatory variable itself.

Furthermore, there is additional information revealed during ongoing processes that cannot be captured easily by such models. An example is concurrency and the effect that new events have on future events. In the online auction context, incoming bids can influence bidders in different ways: A new bid placed in an auction can result in an immediate response by other bidders or can be completely ignored. Bidders also learn from each other: they adopt bidding strategies of other bidders and they learn about an item's value from bids that were placed. Many items sold in online auctions do not have a commonly known value (such as collectibles, antiques, rare pieces of art, etc.), and therefore bidders often try to infer the item's value from other people's bids. In short, while the final price is certainly affected by directly observable phenomena (such as the number of competing bidders), it is also dependent on indirect actions, reactions and interactions among bidders.

### 4.2 Process Dynamics and FDA

Modeling the effects of user interactions with classical regression models is challenging, to say the least. An alternative approach is to capture some of this dynamic information via *evolution* curves and their dynamics. In the auction context this would be the price evolution, which is the progression of bids throughout an auction. The evolution curve and its dynamics can reflect these bidder interactions: High competition in an auction will manifest itself as a steep price curve with increasing dynamics. Price will also increase, albeit at a slower rate, if bidders merely use the new bid to update their own valuation about the product. The price increase will slow down if bidders drop out of the auction due to a newly placed bid or for some other reason. Therefore, the price-evolution curve, and in particular its dynamics, has the ability to capture much of the auction information that would otherwise not be integrated into the model. By price dynamics we mean, for example, the price velocity and acceleration which measure the change in price and the rate at which this change is occurring. The ability to measure dynamics is one of the most noteworthy features of functional data analysis. FDA recovers the price evolution via a smooth curve through the auction's bid history and yields the price dynamics via the derivatives of this curve. Examples of exploring process dynamics via FDA in eCommerce are the price dynamics in eBay auctions (Jank and Shmueli, 2005) and bid dynamics in auctions for modern Indian art by Reddy and Dass (2006). In these two examples the price curves themselves are not very illuminating, but their dynamics reveal interesting patterns and sources of heterogeneity across records.

One can model the relationship between the process evolution (or its dynamics) and other predictors via functional regression analysis. For example, in a few studies of price formation in eBay auctions (Shmueli and Jank, 2006; Bapna, Jank and Shmueli, 2004; Alford and Urimi, 2004) and other online auctions (Reddy and Dass, 2006) a functional regression model was fit to price-evolution curves from eBay auctions (the response) with static predictors (such as the seller rating) and functional predictors (such as the cumulative number of bids). One interesting finding is that the impact of the opening bid on the current price starts high, and slowly decreases as the auction progresses. This reflects the shift in information about the item's value due to bidding: At first there is not much information available and so the opening bid gives a sense of the item's value. But as the auction progresses new bids add more information about the value of the item, thereby reducing the usefulness of the information contained in the opening bid.

One of the challenges in functional regression analysis is the interpretation of the results. Instead of scalar coefficient estimates, we obtain estimated coefficient curves. Plotting these curves means that the $x$-axis is time, and not the ordinary predictor value. For example, Figure 4 shows the estimated coefficient (and a 95% confidence band) for a regression model with a functional response. In the top panel the response is the price evolution. The coefficient is positive throughout the auction, signifying that the current price is positively associated with the opening bid throughout the auction. However, this relationship decreases in magnitude as the auction proceeds. This is reasonable, because bidders gain more and more information as the auction proceeds and therefore derive less utility from the value of the opening price. The middle and bottom panels in Figure 4 describe another useful information source: the relationship between the opening price and the price dynamics. If we are interested in relationships between various independent variables and the process dynamics, we can use the derivative curves as the functional response. In this example we set the price velocity (middle) and price acceleration (bottom) as the responses. We see that the price acceleration is positively associated with the opening bid at the auction start, but then this relationship loses

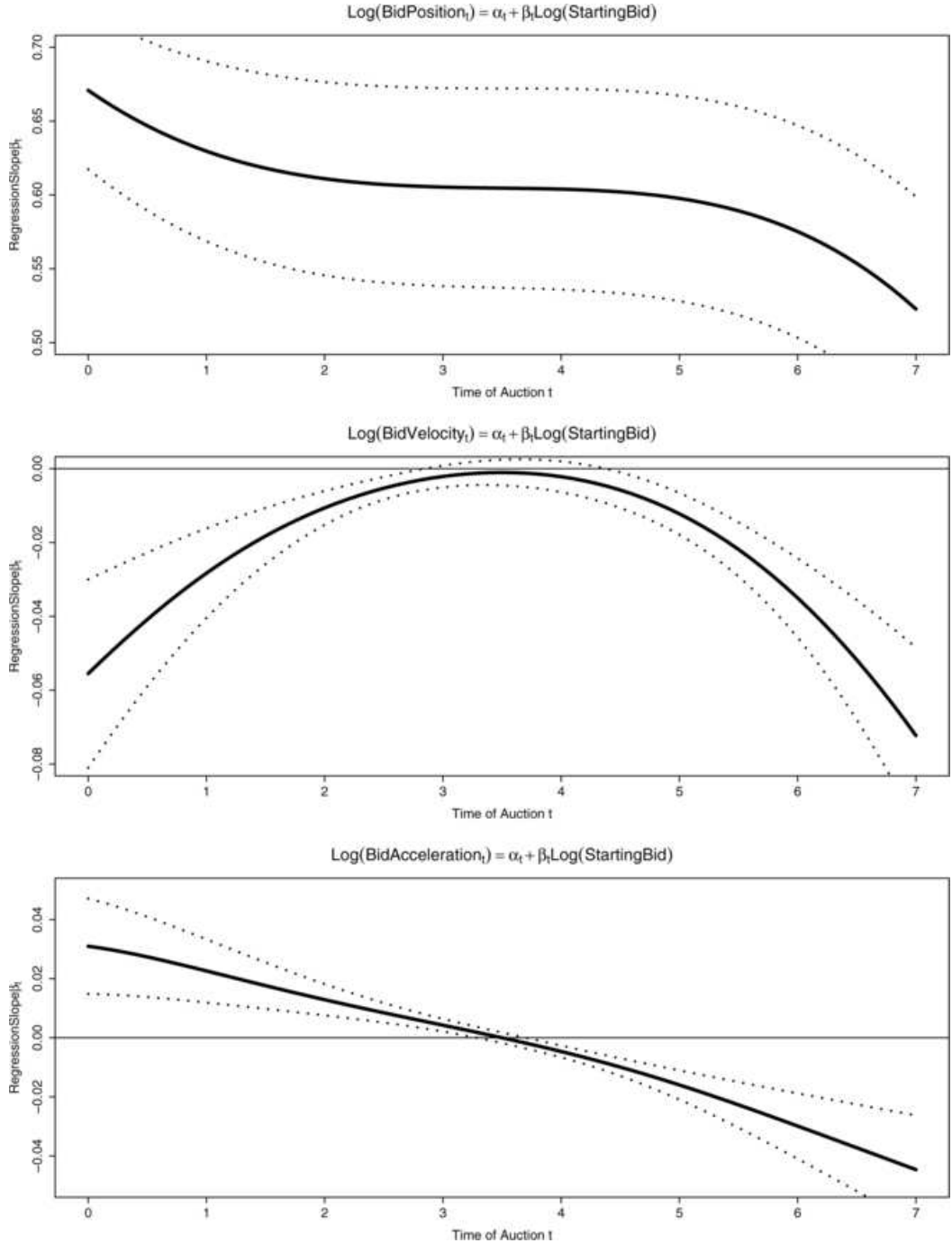

Fig. 4. *Estimated coefficient for opening price, in three regression models with a functional response: price evolution* (top), *price velocity* (middle) *and price acceleration* (bottom).

momentum and even becomes negative as the auction comes to a close.

The capability of studying process dynamics could have a huge impact on eCommerce research. Though the concepts of dynamics are well grounded in physics and engineering, their exact economic impact re-quires more thought. However, there is an opportunity to create new economic measures with the help of FDA. For instance, we can develop concepts such as "auction energy" using the definition of kinetic energy (energy = mass × velocity$^2$/2) to arrive at

auction energy

(1)
$$= (\text{current price}) \times (\text{price velocity})^2/2.$$

A major challenge is to find a theoretical foundation in economics of such concepts. This is only one example where collaboration could have a large impact on the field.

### 4.3 Other Functional Models for eCommerce

Another level of flexibility, but also complexity, is to incorporate interaction terms into functional regression models. Since interactions (in ordinary linear regression models) are widely used in eCommerce studies, it could be useful to measure similar effects in functional objects. The literature on interactions in functional linear regression models appears to be scant, although this seems like an important extension.

Another important direction for modeling the dynamic nature of web content in general, and eCommerce in particular, is the use of differential equations. The use of differential equations in the functional literature is still in its infancy. Ramsay (2000a) gives an introduction to the use of differential equations in statistics and several examples of functional estimation problems such as simultaneous estimation of a regression model and residual density, monotone smoothing, specification of a link function, differential equation models of data, and smoothing over complicated multidimensional domains. Ramsay calls this "principal differential analysis" (PDA) because of the similarities that it shares with principal components analysis.

PDA is a natural formalization of the exploration of curve dynamics. Through a differential equation we can relate, for instance, the price during an auction to its rate of increase and acceleration. It is possible that such relations exist in the dynamic, ever-changing eCommerce world. These relationships need to be more formally integrated with economic theory to create a solid foundation for the empirical findings that have been observed. For example, PDA was used to study the relationship between price curves in online auctions and their derivatives by Jank and Shmueli (2005) and Wang (2005). The main finding is that relationships between the price curve and its acceleration are present in some types of auctions, but not in others, suggesting that dynamics can vary broadly in eCommerce processes.

## 5. FUTURE TRENDS IN FUNCTIONAL MODELING

In the previous sections we have shown multiple facets of FDA that make it a natural approach in eCommerce empirical research. Unlike currently used static models, FDA can capture processes and dynamic information which are inherent in the eCommerce environment. In the following we describe a few important areas that are still undeveloped both in the eCommerce research world and in the FDA domain, and in our opinion have the potential to make a contribution to both.

The first area is related to concurrency of events. In almost every eCommerce study, the events of interest occur concurrently or have at least some overlap. This means that there is a dependence structure between records, with some events influencing others. The most obvious example is the stock market with stock prices influencing each other. Some eCommerce examples are concurrent auctions on eBay for the same item or even for competing items, and prices of a certain book at different online vendors (and perhaps even brick-and-mortar stores) over time. Although researchers acknowledge such relationships, nearly all studies make the simplifying assumption of independent observations. Ignoring the effects of concurrency can lead to invalid results. A first step is therefore to find ways to evaluate the degree of concurrency and its effect on the measure of interest. Shmueli and Jank (2005) introduce and evaluate several data displays for exploring the effect of concurrency in online auctions on final price. Hyde, Jank and Shmueli (2006) expands this work to visualize concurrency of the functional objects, using curves to represent price evolution and its dynamics. In addition to data displays, there is a need for defining measures of concurrency, and finally, for developing models that can incorporate and account for relationships between processes that are represented as functional objects (see, e.g., Jank and Shmueli, 2006).

Another important enhancement to FDA that would greatly benefit eCommerce research is the incorporation of change into the functional objects over time. As pointed out earlier, eCommerce experiences constant change. Frequent technological advancement, new website formats, changes in the global economy, etc., can have a large influence on what we observe in the eCommerce world. If we use functional objects to represent observations which are themselves longitudinal, we need ways to incorporate an additional temporal dimension that compares functional objects over time.

FDA research focuses more and more on $p$-dimensional functional objects (e.g., Yushkevich et al., 2001). In many eCommerce applications such representations could be very useful. One example is competition in online auctions, where each auction is represented by its price curve coupled with the cumulative number of bidders, thus yielding a bivariate functional representation (Wang and Wu, 2004). There is also related work on *symbolic data analysis* (SDA) (see Bock and Diday, 2000), which provides tools for managing complex, aggregated, relational and higher-level data described by multivalued variables. This could be a new successful wedding with FDA methods.

Finally, in many cases, and especially in economics, the objects of interest are individual users. Economists are typically interested in how individuals strategize and react to others. The problem with formulating functional objects that represent individuals is sparsity of data. That is, individuals typically do not leave many traces during one eCommerce transaction. For example, in eBay auctions if we treat an individual bidder as our object of interest, then bidders will leave very sparse data (1–2 bids per bidder is the norm). In a similar setting, James, Hastie and Sugar (2000) and James and Sugar (2003) use a semiparametric setting where information from data aggregated across individuals is used to supplement the information at the individual level. An alternative model would pool information from previous records that the individual was involved in and use it to supplement the current record. Approaches that enable the functional representation of sparse data can prove very useful in tying economic theories to empirical results. This would further strengthen the eCommerce research area.

## 6. CONCLUSIONS

The emerging field of empirical eCommerce research is growing fast with many data-related challenges. In light of the special data structures and the types of research questions of interest, we believe that functional data analysis can play a major role in this field. On the one hand, this requires more involvement by statisticians to further explore statistical issues involved and to develop functional methods and models that are called for in these applications. On the other hand, collaborative work has proven to be extremely fruitful for the multiple disciplines involved. In that respect, more outreach should be done to make these tools more popular. Wider adoption of functional tools by nonstatisticians requires software accessibility. Currently FDA packages exist for Matlab, S-PLUS and R (ego.psych.mcgill.ca/misc/fda/software.html). Specialized programs for particular applications are anticipated to grow, and we encourage researchers to make such code and data freely available. In addition, making sample datasets freely available will make this field more accessible and attractive to statisticians. Our website (www.smith.umd.edu/ceme/statistics/) contains some eBay data and auction-specific FDA code.

Another important front in further developing this exciting new interdisciplinary field is the involvement and training of graduate students. This includes educating statistics students about both the eCommerce domain and FDA. From our own experience through Interactive Research Teams, we found this to be a very exciting ground for advancing statistical research.

## ACKNOWLEDGMENTS

We thank Professor Steve Marron from UNC for fruitful conversations, Professor Jim Ramsay from McGill University for continuous advice and support with FDA software, and the three referees for their constructive comments.